\documentclass[12pt,a4paper]{amsart}
\usepackage{graphicx}
\usepackage{amscd}  
\newtheorem{theorem}{Th\'{e}or\`{e}me}
\theoremstyle{plain}

\newtheorem{corollary}{Corollaire}

\newtheorem{definition}{Definition}

\newtheorem{assume}{Hypoth\`{e}se}
\newtheorem{remerciements}{Remerciements}

\newtheorem{proposition}{Proposition}
\newtheorem{remark}{Remarque}

\numberwithin{equation}{section}
\numberwithin{theorem}{section}
\numberwithin{proposition}{section}
\numberwithin{lemma}{section}

\newcommand{\supp}{\operatorname{supp}}

\newcommand{\tr}{\operatorname{tr}}

\renewcommand{\Re}{\mathfrak{R}}
\renewcommand{\Im}{\mathfrak{I}}

\numberwithin{equation}{section}
\begin{document}
\title{Calcul de Wick en dimension infinie}
\author{Bernard Lascar}

\address
{Bernard Lascar. Institut de Math\'ematiques.  Analyse Alg\'ebrique. Tour 15-25 5. Universit\'e Pierre et Marie Curie. 2 Place Jussieu 75005. Paris. France.}
\footnote{avec une annexe de R.~Lascar {\tiny Richard.Lascar@unice.fr}}
\begin{abstract}On prouve ici des r\'esultat de bornes inf\'erieures
  \`a la Melin-H\"ormander en grandes dimension par la m\'ethode du
  calcul de Wick.
\end{abstract}
\maketitle
Soit $m\in\mathbb{R}$ et $k\in\mathbb{N}$, $k\geq 1$, et $\Lambda\geq
1$ le param\`etre spectral; on note 
\begin{definition}
\label{def1}
Soit $\mathcal{N}^{m,2k}(\mathbb{R}^n,\Sigma)$ la classe des
op\'erateurs pseudo-diff\'erentiels classiques de symboles
$p(x,\xi,\Lambda)$ qui ont la propri\'et\'e :
\begin{equation}
\label{0}\begin{split}
  &P(x,\frac{D_{x}}{\Lambda},\Lambda)u(x)=\\&\iint
  p(\frac{x+y}{2},\xi,\Lambda)\exp{(i\Lambda(x-y)\xi)}u(y,\Lambda)dy\left(\frac{\Lambda
      d\xi}{2\pi}\right)^n,\\&\hskip 5cm \text{o\`u $u\in
    \mathcal{S}(\mathbb{R}^n)$}\end{split}
\end{equation}
\begin{equation}
\label{1}
p(x,\xi,\Lambda)\sim \Lambda^{m}p_m(x,\xi)+\ldots+\Lambda^{m-j}p_{m-j}(x,\xi)+\ldots
\end{equation}
o\`u $p_{m-j}\in S(1,\Gamma)$, avec $\Gamma=|dX|^2$ et
$p_{m-j}(X)=\mathcal{O}(d_{\Sigma}^{(2k-2j)_{+}})(X)$. Avec
$d_{\Sigma}$ la distance euclidienne \`a une sous-vari\'et\'e lisse
symplectique de $\mathbb{R}^n\oplus\mathbb{R}^n$, de codimension $2d$
et donc $d\leq n$.
\end{definition}
Une condition comme $(\ref{1})$ est invariante par transformation
canonique et donc on peut supposer que $$\Sigma=\{(x,\xi);
x''=\xi''=0\},\ (x'',\xi'')\in \mathbb{R}^d\oplus\mathbb{R}^d.$$ On
change par souci de commodit\'e les notations et $X"$ devient $X\in
\mathbb{R}^d\oplus\mathbb{R}^d$ tandis que la variable le long de
$\Sigma$, $X_n\in\mathbb{R}^{n-d}\oplus\mathbb{R}^{n-d}$ est not\'ee
$X_n$. 

On suppose :
\begin{assume}
\label{H1}
\begin{itemize}
\item[i)]  $P^*=P$, $p_{m}\simeq d_{\Sigma}^{2k}.$
\item[ii)] En tout point $\rho\in\Sigma$, le symbole de Wick du
  localis\'e $P_{\rho}$ est un symbole uniform\'ement avec $(\rho,
  \Lambda,d)$ elliptique dans $S(\Lambda^{m-k}(1+|z])^{2k},G)$ o\`u
  $G=\frac{|dz|^{2}}{(1+|z|)^{2}}$
\end{itemize}
\end{assume}

\begin{theorem}
\label{th}
Sous les hypoth\`eses $(H_1)$, il existe deux constantes $C_0>0$ et $C_1>0$ telles que :
\begin{multline}
\label{2}
(Pu,u)\geq -C_0 \Lambda^{m-k-1}\| u\|^2,\ \text{pour $u\in
  \mathcal{S}(\mathbb{R}^n)$}\\ \text{uniform\'ement pour
  $d\Lambda^{-1}\geq C_1$ et $(n-d)\leq C_1^{-1}$,}\\ \text{ou
  bien-s\^ur si $d$ est fixe et $\Lambda\rightarrow\infty$.}
\end{multline}
\end{theorem}
\begin{remark}
Quand $k=1$ il s'agit du th\'eor\`eme \cite{Lascar}.
\end{remark}

On effectue comme dans \cite{Lascars} une dilatation pr\`es des
caract\'eristiques,
\begin{multline}
\label{3}
u(x,x_n)\rightarrow
v(y,y_n,\Lambda)=u(\Lambda^{-1/2}y,y_n)\Lambda^{-d/4},\ \text{dans le domaine} \\
 \text{$E_{\varepsilon}=\{(Y,Y_n); |Y|\leq \varepsilon \Lambda^{1/2}, Y_n\in V(\rho_0)\}$, $(Y,Y_n)\in \mathbb{R}^n\oplus\mathbb{R}^n$}
\end{multline}
En effet, lorsqu'on change $x$ en $\Lambda^{1/2}x$, $\xi$ est chang\'e
en $\Lambda^{-1/2}\xi$ mais comme on quantifie le nouvel op\'erateur
pseudo-diff\'erentiel en Weyl-$(1,\Lambda)$, c'est que $\xi$ devient
$\Lambda^{1/2}\xi$, on retrouve donc une bonne homog\'en\'eit\'e, et
c'est le symbole de Weyl qui est transform\'e exactement.

\section{L'argument de Melin-H\"ormander en dimension $n$.}

\begin{multline}
\label{4}
(Qv)(y,y_n,\Lambda)=\iint q(\frac{x+y}{2},\eta,\Lambda)v(y,\Lambda) \\
\exp{[i(x-y)\eta+\Lambda (x_{n}-y_{n})\eta_{n}]}\left(
  \frac{\Lambda}{2\pi}\right)^{n-d}dy_n
d\eta_n\left(\frac{dyd\eta}{2\pi}\right)^d \ \\ \text{o\`u
  $q(Y,Y_n,\Lambda)=p(\Lambda^{-1/2}Y,Y_n,\Lambda)$.}
\end{multline}
On rappelle que L. Boutet de Monvel \cite{Bo} a introduit pour $m, k\in\mathbb{R}$ la classe :
\begin{multline}
\label{5}
\tilde{S}^{m,k}(\mathbb{R}^n,\Sigma)=\\ \{\tilde{f}\in C^{\infty}(\mathbb{R}^n); |D_X^pD_{X_{n}}^q\tilde{f}X,X_n)| \leq \Lambda^{m}d_{\Sigma}^{k-p}(X,X_n)\}\  \\ \text{avec $d_{\Sigma}(X,X_n)=\Lambda^{-1/2}+|X|$}.
\end{multline}
La m\'etrique correspondante \cite{Ho} est 
\begin{equation}
\label{6}
\tilde{g}_{X}=\frac{|dX|^{2}}{d_{\Sigma}^{2}(X)}+|dX_n|^2
\end{equation}
Le poids est $\tilde{m}(X)=\Lambda^{m}d_{\Sigma}^k(X)$. Donc la dilatation produite en (\ref{3}) produit la classe $S(\Lambda^{m-k/2}d(Y)^{k},g)$ pour $f(Y,Y_n)=\tilde{f}(\Lambda^{-1/2}Y,Y_n)$ avec
\begin{equation}
\label{7}
d_a(Y)=(a^2+|Y|^2)^{1/2},\ \text{et $g_{Y}=\frac{|dY|^{2}}{d_a^{2}(Y)}+|dY_n|^2$, avec $a=1$.}
\end{equation}

On voit facilement que :
\begin{equation}
\label{7.1}
|d_a(X+Y)-d_a(X)|\leq |Y|, \ \text{$|\frac{d_{a}(X+Y)}{d_{a}(X)}-1|\leq g_{a,X}(Y)^{1/2}$}
\end{equation}
Donc les m\'etriques $g_a$ sont lentes d\`es que $a>0$.

Puis :
\begin{equation}
\label{7.2}
\frac{1}{d_{a}^{2}(X)}\leq \frac{1}{d_{a}^{2}(X+Y)}2(1+g_{a,X}(Y))
\end{equation} 

La temp\'erance r\'esulte donc de la lenteur et de $h_{a}\leq 1$ qui est assur\'ee d\`es que $a\geq 1$.
 \begin{equation}
 \label{7.3}
 h_a(Y)=\max({d_{a}^{-2},\Lambda^{-1})}
 \end{equation}
 
\begin{equation}
\label{8}
\begin{split}
\omega&=\sigma_Y+\Lambda^{-2}\sigma_{Y_{n}},\text{ et donc}\\
g_a^{\omega} &=d_a^2|dY|^2+\Lambda^{-2}|dY_n|^2\geq \min{ (d_{a}^{-4},\Lambda^{-2})}g_a
\end{split}
\end{equation}
Il faut donc que $d_a(Y)\geq 1$, ce qui donne bien $a\geq 1$.

Une autre remarque est que si $\tilde{f}(X,X_n)\in S(1,\Gamma)$ alors
$f(Y,Y_n)\in S(1,\Lambda^{-1}|dY|^2+|dY_n|^2)$, il faut donc imposer
aussi la condition $\Lambda^{-1}\leq C_1d_{a}^{-2}$, soit aussi , ce
qui impose de se limiter \`a une zone comme $E_C$ et de prendre aussi
$1\leq a\leq C\Lambda^{1/2}$. Alors $\tilde{f}\in
S(1,\Gamma)\Rightarrow f\in S(1,g)$, au moins dans $E_C$. Et donc
$h_a\sim d_a^{-2}$.   

Le calcul mixte d de classes
$S(m,g_a)\#S(m'g_b)\subset S(mm',g_{\min{(a,b)}})$ et $h_{a,b}\sim
d_a^{-1}d_b^{-1}$.
 
 Si on part de $\tilde{f}\in
S(\tilde{m},g_{\tilde{b}})$, alors $f\in
S(m,g_{\Lambda^{1/2}\tilde{b}})$, avec
$m(Y)=\tilde{m}(\Lambda^{-1/2}Y)$.
\begin{equation}
\label{8.0}
\tilde{f}\in S(\Lambda^{m}d_{\mu}^{\beta},g_{\mu})\Rightarrow f\in S(\Lambda^{m-\beta/2}d_{\mu\Lambda^{1/2}}^{\beta},g_{\Lambda^{1/2}\mu})
\end{equation}

En particulier si $\tilde{f}\in S^{\alpha}=S(\langle X\rangle^{\alpha},\frac{|dX|^{2}}{\langle X\rangle^2})$, alors $f\in S(d_{\Lambda^{1/2}}^{\alpha}\Lambda^{-\alpha/2},g_{\Lambda^{1/2}})$. Ce qui fait que si $\tilde{f}\in \mathcal{S}(\mathbb{R}^d\oplus\mathbb{R}^d)$, $f\in\bigcap_{\alpha\in \mathbb{R}^{-}}S(d_{\Lambda^{1/2}}^{\alpha}\Lambda^{-\alpha/2},g_{\Lambda^{1/2}})$.

Comme on peut toujours partir d'un probl\`eme microlocalis\'e au sens
usuel pr\`es de $\Sigma$ on sera toujours pour $q(Y,Y_n,\Lambda)$ dans
la classe :
\begin{equation}
\label{8.1}
q(Y,Y_n,\Lambda)\in \bigcap_{\alpha\in \mathbb{R}^{-}}S(d_{\mu\Lambda^{1/2}}^{\alpha}\Lambda^{m-\alpha/2},g_{\mu\Lambda^{1/2}})\overset{\text{def}}{=}\mathcal{T}_{\mu}^m(\mathbb{R}^n,\Sigma)\ \text{avec $\mu=1$}.
\end{equation}

On choisira convenablement $\mu>0$ et une partition de l'unit\'e $(\chi_{\mu})_{\mu\in N}$, $\chi_{\mu}\in S(1,g_{\mu})$, telle que 
\begin{equation}
\label{8.2}
\begin{split}
&\sum_{\mu\in N}\chi_{\mu}^2(X)=\zeta_0(X),\text{ o\`u }\zeta_0\in
C^{\infty}_0(\mathbb{R}^n\oplus\mathbb{R}^n), \\ & \zeta_0\equiv 1 \text{
  pr\`es de } (0,\rho_0) \text{ et est support\'ee pr\`es de ce m\^eme point.}
\end{split}
\end{equation}

On aura donc remplac\'e $q(Y,Y_n,\Lambda)$ par $q_{\nu}=\chi_{\Lambda, \mu}q$ avec $\chi_{\Lambda,\mu}(Y,Y_n)=\chi_{\mu}^2(\Lambda^{-1/2}Y,Y_n)$ o\`u les $\chi_{\mu}$ sont donn\'ees par (\ref{8.2}). Le nouveau grand param\`etre de deuxi\`eme microlocalisation est de la forme $\nu=\Lambda^{1/2}\mu$, avec $\nu$ grand et $\mu=\Lambda^{-1/2}\nu$ petit.

Si $\varepsilon\leq 1$, on peut former le d\'eveloppement de Taylor \`a l'ordre $N$ 
\begin{multline}
\label{9}
f(Y,Y_n)=\sum_{0\leq j\leq N-1}\frac{1}{j!}D_Y^jf(0,Y_n)\ Y^j+R_N((f)Y,Y_n)\  \\ \text{avec $R_N(f)(Y,Y_n)=\int_0^1\frac{(1-\tau)^{N-1}}{(N-1)!}D_Y^Nf(\tau Y,Y_n)Y^Nd\tau$.}
\end{multline}
Comme $Y\in S(d_a,g_a)\forall a$; et que si $\tilde{f}\in S(1,\Gamma)$ alors 
\begin{equation}
\label{10}
\forall C>0 \ \text{dans $E_C$, $R_N(f)(Y,Y_n)\in S(\Lambda^{-N/2}d_a^N,g_a)$.}
\end{equation}
Si maintenant $\tilde{f}\in \mathcal{N}^{m,2k}$, alors $R_N(f)(Y,Y_n)\in S(\Lambda^{m-N/2}d_a^N,g_a)$ et $D_Y^jf(0,Y_n)=0$ pour $j\leq 2k-1$, donc 
\begin{equation}
\label{11}
f(Y,Y_n)=\sum_{2k\leq j\leq N-1}R_j(f)(0,Y_n)+R_N(f)(Y,Y_n)
\end{equation}
Les conditions $d_a\leq C\Lambda^{1/2}$ entra{\^\i}nent alors que $f(Y,Y_n)\in S(\Lambda^{m-k/2}d_a^k,g_a)$ et $R_N(f)(Y,Y_n)\in S(\Lambda^{m-N/2}d_a^{N},g_a)$ pour tout $N$ et $R_j(f)(0,Y_n)=0$ pour $j<2k$.

\begin{multline}
\label{12}
q(Y,Y_n)=\sum_{0\leq l\leq k}\left(\sum_{2(k-l)\leq j\leq N_{l}-1}R_{j}(q_{m-l})(0,Y_n)+R_{N_{l}}(q_{m-l})(Y,Y_n)\right)+\\ q_{m-k-1}+\ldots+q_{m-k-M}+\ldots
\end{multline}
Si on choisit $N_{l}=2(k-l)+1$ dans (\ref{12}), on obtient :
On a donc 
\begin{proposition}
\label{prop1}
Si $P\in\mathcal{N}^{m,2k}$ alors 
\begin{multline}
\label{13}
q(Y,Y_n)=\sum_{0\leq l\leq k}\left(R_{2(k-l)}(q_{m-l})(0,Y_n)+R_{2(k-l)+1}(q_{m-l})(Y,Y_n)\right)\\ +q_{m-k-1}+\ldots+q_{m-k-M}+\ldots
\end{multline}
et donc 
\begin{multline}
\label{14}
q(Y,Y_n)=\\\sum_{0\leq l\leq k}\left(R_{2(k-l)}(q_{m-l})(0,Y_n)+R_{2(k-l)+1}(q_{m-l})(Y,Y_n)\right)+S(\Lambda^{m-k-1/2}d_a^{2k+1},g_a)\ \\
\text{o\`u $a$ est arbitraire pourvu que $1\leq a\leq
  \Lambda^{1/2}$.}
\end{multline}
\end{proposition}
Le Corollaire de (\ref{14}) \`a l'aide des arguments de seconde microlocalisation \'evoqu\'es plus hauts conduisent au th\'eor\`eme optimal \cite{Lascars} suivant:
\begin{theorem}
Sous les hypoth\`eses :
\label{th1}
\begin{assume}
\label{H2}
\begin{itemize}
\item[i)]  $P^*=P$, $p_{m}\simeq d_{\Sigma}^{2k}.$
\item[ii)] En tout point $\rho\in\Sigma$, l'op\'erateur localis\'e $P_{\rho}\geq 0$ sur $L^2$.
\end{itemize}
\end{assume}
on d\'eduit :  pour tout $n\in\mathbb{N}$, il existe $C_n>0$ tel que 

\begin{equation}
\label{15}
(Pu,u)\geq -C_{n}\Lambda^{m-k-1/2}\|u\|^2\ \text{$\forall u\in\mathcal{S}(\mathbb{R}^n)$}.
\end{equation}
\end{theorem}

\section{Espace de Fock et calcul de Toeplitz}

On travaille directement sur l'espace de Fock commutatif d'un espace
de Hilbert complexe $H$ de dimension quelconque, dont on supposera
qu'il est le complexifi\'e d'un espace de Hilbert r\'eel $H_1$
s\'eparable.

On note $|z|$ la norme de $H$ et pour $k\in\mathbb{N}$, l'espace
$W_k=\otimes_k H$, et $V_k$ le sous-espace engendr\'e par les tenseurs
sym\'etriques, l'op\'erateur de sym\'etrisation total op\`ere bien de
$W_k\otimes W_l$ dans $V_{k+l}$ est a bien pour norme $1$.

Si $n=\dim H$ il est clair que $\dim
V_k=\frac{(n+k)!}{n!k!}=\mathcal{O}(N^{k})$ quand $n\rightarrow
\infty$. Donc plus la dimension $N$ augmente plus $V_k$ verra sa
dimension augment\'ee avec $k$.

L'espace de Fock $F(H)=\oplus_k V_k$ avec la convention de somme hilbertienne :
\begin{multline}
\label{1}
\Phi=\sum_k \Phi_k, \ \Phi_k\in V_k  \ \Phi_k(z)=\sum_{|\alpha|=k}\Phi_{\alpha}\frac{z^{\alpha}}{\alpha !^{1/2}}\\ \text{et donc }\|\Phi \|^2=\sum_{\alpha}|\Phi_{\alpha}|^2=\int_H |\Phi (z)|^2 d\nu' (z)
\end{multline}

Dans $(\ref{1})$, $\nu'$ est la mesure gaussienne sur $H$ de variance $1/2$, $\nu ' =\exp(-|z|^2)\frac{L(dz)}{\pi^n}$ quand $\dim H=n$.

Donc $F(H)=L^2(H,\nu')\cap\mathcal{E}$, o\`u $\mathcal{E}$ est l'espace des fonctions enti\`eres, avec juste une convention \`a expliquer quand $\dim H=\infty$ car $H$ est de mesure nulle dans $L^2(H,\nu ')$ et donc ne peut avoir v\'eritablement de fonctions continues et enti\`eres les simples polyn\^omes comme $|z|^{2k}$ \'etant presque partout infinis d\`es que $k\geq 1$ ce qui est la source v\'eritable des infinis de renormalizations.

Il y a bien un projecteur orthogonal de $L^2(H,\nu')$ sur $F(H)$ not\'e par $\pi$ donn\'e par la formule de la TFBI \cite{Sj} : 
 
\begin{multline}
\label{2}
\pi f (z)=\int_H f(w,\overline w) \exp (z \overline w)d\nu' (w);\\ \pi (f)=f \text{ si et seulement si  $f$ est holomorphe} 
\end{multline}
$e_z=\exp( z \overline w)$ est bien pour $z\in H$ est dans $L^2_w$ et donc l\`a il n'y a plus d'ambiguit\'e de diff\'erentiabilit\'e. $e_z$ est l'\'etat coh\'erent en $z$.

L\`a c'est purement $L^2$ car une difficult\'e connue (calcul de Malliavin par exemple) est que $\pi$ n'a plus de raisons d'op\'erer dans $L^p$. $L^p$ pour $p\neq 2$ est donc une autre th\'eorie.

On peut d\'efinir des op\'erateurs non triviaux dans $F$ par la formule :

\begin{equation}
\label{3}
A(\Phi)=\pi (a^{anti-wick}\Phi) \ \text{o\`u $a^{anti-wick}\in L^{\infty}_{\nu^{\prime}}$}
\end{equation}

$\|A\| \leq \| a^{anti-wick}\|_{L^{\infty}_{\nu^{\prime}}}$ avec une in\'egalit\'e stricte en g\'en\'eral car $F$ n'engendre pas $L^2$. 

Nous dirons que $a^{anti-wick}$ est l'anti-symbole de Wick de $A$, il est bien s\^ur uniquement d\'efini par $(\ref{3})$ on note alors $A=(a^{anti-wick})^{w}$ et on d\'efinit le symbole de Wick 

$a^{wick}(z,\overline z)=\tilde{ \pi}(a^{anti-wick})(z,\overline z)$. 

$\tilde{\pi}$ est encore une TFBI \`a deux variables cette fois que l'on restreint \`a la diagonale $\{(w,z)\in H\times H;w=\overline z\}$, $a^{wick}(z,\overline z)$ est lui bi-holomorphe et dans $L^2(H,\nu ')\otimes L^2(H,\nu ')\cap\mathcal{E}_1$, cette fois $\mathcal{E}_1$ est l'espace des fonctions bi-holomorphes en $(z,\overline z)$ et on a aussi l'autre \'egalit\'e : 

\begin{equation}
\label{4}
\int_{H\times H}|a^{wick}(z,\overline w)|^2 d\nu' (z)d\nu' (w)=\int_H |a^{anti-wick}(w,\overline w)|^2 d\nu ' (w).
\end{equation}

Soit $\mathcal{H}_{k,l}$ le polyn\^ome d'Hermite obtenu par la formule : 
\begin{equation}\label{5}
 \exp(-|z+w|^2)=\sum_{k,l\in \mathbb{N}}(-1)^{k+l}\mathcal{H}_{k,l}(w,\overline w)\frac{\overline{z}^{k}z^{l}}{k!l!}\exp(-|w|^2)
\end{equation}
Les polyn\^omes d'Hermite obtenus par (\ref{5}) sont deux \`a deux orthogonaux dans $L^2(H,\nu')$ et les $H_{k,l}=\mathcal{H}_{k,l}(k!l!)^{-1/2}$ sont la base orthonorm\'ee de $L^2(H,\nu')$ obtenue par orthonormalisation des $z^k\overline{z}^l$. Il r\'esulte alors de $(\ref{5})$ que 
\begin{multline}
\label{6}
\tilde{\pi}(f)(z,\overline z)=\exp(-|z|^2)\int f(w,\overline
w)\exp(z\overline w+\overline z w)d\nu'(w)= \\
\sum_{k,l}(-1)^{k+l}\int f(w,\overline
w)\overline{\mathcal{H}_{k,l}}(w,\overline
w)\frac{\overline{z}^{l}z^{k}}{k!l!}d\nu'(w)\\ = \sum_{k,l}(-1)^{k+l}
f_{k,l}\frac{z^{k}\overline{z}^l}{(k!l!)^{1/2}}\ \text{o\`u
  $f_{k,l}=(f,H_{k,l})$}
\end{multline}

Ce qui prouve $(\ref{4})$ et aussi que :
\begin{equation}
\label{7}
\tilde{\pi}(f)(z,\overline z)=\exp(\partial_{\overline{z}}\partial_z)f(z,\overline{z})=\sum_{\theta\in \mathbb{N}^{\mathbb{N}}}\frac{1}{\theta !}\partial_z^{\theta}\partial_{\overline z}^{\theta}f(z,\overline z)
\end{equation}
La formule (\ref{7}) s'inverse par :
\begin{equation}
\label{8}
f(z,\overline z)=\exp(-\partial_{\overline{z}}\partial_z)f(z,\overline{z})=\sum_{\theta\in \mathbb{N}^{\mathbb{N}}}(-1)^{|\theta|}\frac{1}{\theta !}\partial_z^{\theta}\partial_{\overline z}^{\theta}\tilde{\pi}(f)(z,\overline z)
\end{equation}
Les formules $(\ref{6})$ et $(\ref{7})$, montre que si l'op\'eration $f\rightarrow \tilde{\pi}(f)$ ne pose pas de probl\`eme, l'op\'eration inverse ne se produit pas dans les distributions temp\'er\'ees.

Le calcul du compos\'e $C$ de $B$ et $A$; o\`u $B$ et $A$ sont des
op\'erateurs ayant pour anti-symboles de Wick respectivement les
fonctions $b^{anti-wick}$ et $a^{anti-wick}$, au moins formellement
$C$ est lui donn\'e par la formule :
\begin{equation}
\label{9}
c^{anti-wick}(z,\overline z)=\sum_{j\in \mathbb{N}}\frac{(-1)^{j}}{j!}\tr(D_{z}^jb^{anti-wick})\otimes (D_{\overline z}^{j}a^{anti-wick})(z,\overline z)).
\end{equation}

La formule $(\ref{9})$ est exacte quand $a^{anti-wick}$ ou $b^{anti-wick}$ est un polyn\^ome.

\section{Op\'erateurs de cr\'eation et d'annihilation dans l'espace de
  Fock.}

Avant de regarder des op\'erateurs de Toeplitz abstraits on consid\`ere les annihilateurs $a_{k}=D_z^{k}$ et les cr\'eateurs $a^{*}_j=z^{j}$.
On d\'efinit un espace de Sobolev
\begin{definition}
\label{defS}
Soit $s\in\mathbb{R}$, l'espace de Sobolev $F^{s}(H)$ est l'ensemble
des $\Phi (z)\in \mathcal{E}$ telles que :
\begin{equation}
\label{2.1}\begin{split}
&\Phi=\sum_k \Phi_k, \ \Phi_k\in V_k,
\Phi_k(z)=\sum_{|\alpha|=k}\Phi_{\alpha}\frac{z^{\alpha}}{\alpha
  !^{1/2}}\\ &\text{et }  \| \Phi \|^2_s=\sum_k |\Phi_k|^2 (1+k)^s <\infty
\end{split}
\end{equation}
\end{definition}

On constate que 
\begin{proposition}
\label{prop2.1}
$a_{k}\in \mathcal{L}(F^s,F^{s-k}\otimes V_k)$ et $\| a_{k}\| \leq \sup_{ q\geq k} (1+q-k)^{s-k}(1+q)^{-s}\frac{q!}{(q-k)!k!}$.
\end{proposition}
Donc par passage \`a l'adjoint :
\begin{proposition}
\label{prop2.2}
$a_{k}^{*}\in \mathcal{L}(F^s\otimes V_k,F^{s-k})$ et $\| a_{k}^{*}\| \leq \sup_{ q\geq k} (1+q-k)^{s-k}(1+q)^{-s}\frac{q!}{(q-k)!k!}$.
\end{proposition}

On peut construire des oscillateurs harmoniques d'ordre $m$ par 
\begin{equation}\label{2.2}\begin{split}
&A=\sum_{k+l\leq m} a_{k}^{*}a_{k,l}a_{l},\ a_{k,l}\in\mathcal{L}(V_l,
V_k) \\&\text{et ainsi }A\in \mathcal{L}(F^s,F^{s-m}) \text{ pour tout }s\in\mathbb{R}.
\end{split}
\end{equation}

Il est facile de voir que l'on constitue ainsi une alg\`ebre d'op\'erateurs diff\'erentiels \`a coefficients polynomiaux et que les symboles principaux 
\begin{equation}
a_m(z,\overline z)=\sum_{k+l=m}\tr z^{\otimes k}\otimes a_{k,l}\overline{z}^{\otimes l}
\end{equation}
se multiplient, et que $A^*$ a pour anti symbole de Wick la fonction $\overline {a^{anti-wick}}(z,\overline z)$.

Quand $\dim H=n$, l'anti symbole de Wick de $A(z,D_z)$ est donn\'e par la formule :
\begin{equation}
\label{2.3}
a^{anti-wick}(z,\overline z)=\sum_{k+l\leq m, \tau\leq \min(k,l)}\frac{(-1)^{\tau}}{\tau !}\tr_{\tau}D_z^{\tau}D_{\overline z}^{\tau}a^{wick}(z,\overline z).
\end{equation}

Il r\'esulte de $(\ref{2.2})=$ que pour $0\leq p\leq k$ et $0\leq q\leq l$ 
\begin{equation}
\label{2.4}
\|D_z^p D_{\overline z}^q(a_{k,l}(z^k,\overline z^l) )\|_{\mathcal{L}(V_{p},V_{q}^{*})}\leq C |z|^{k+l-p-q}
\end{equation}
$D_z^{\tau}D_{\overline z}^{\tau}A(z,\overline z)\in {\mathcal{L}(V_{\tau},V_{\tau}^{*})}$, la dimension de $V_{\tau}$ \'etant $\dim V_{\tau}=\frac{(n+\tau)!}{n!\tau !}=\mathcal{O}(n^{\tau})$, 
\begin{equation}
\label{2.4}
\tr D_z^{\tau}D_{\overline z}^{\tau}A(z,\overline z)=\mathcal{O}(|z|^{k+l-2\tau}n^{\tau}).
\end{equation}
\begin{remark}
\label{remark1}
Quand $k=1$ et $l=1$, c'est l'oscillateur harmonique; il y a deux
termes $(az,\overline z)$, $a\in \mathcal{L}(H)$ se comporte comme
$|z|^2$ quand $a$ est inversible positif, et a donc une norme
$L^1_{\text{loc}}$ \'egale \`a $n$, le deuxi\`eme terme est $\tr
a\simeq n$; ils ne se compensent pas, m\^eme en moyenne quadratique
puisque la norme $L^2$ vaut :
\begin{equation}\label{2.5}
\int_H |(az,\overline z)-\tr a)|^2d\nu' (z)=\int_H|\sum_{j=1}^n\lambda_j H_{(j),(j)}|^2d\nu'(z)=\sum_{j=1}^n\lambda_j^2\geq n \| a\|^2
\end{equation}
\end{remark}

\section{Calcul symbolique de Toeplitz}
Soit $\theta\in \mathbb{C}$ soit $N_{\theta}$ l'op\'erateur
d'anti-Wick $ n_{\theta}=(|z|^{2}+n)^{\frac{\theta}{2}}\in S(
(n^{1/2}+|z|)^{\Re \theta},G))$.

On introduit donc la translation symplectique unitaire, $\tau_{Y}$ est
l'op\'erateur :
\begin{equation}
\tau_Y=(a_Y)^{anti-wik}\ \text{o\`u $a_Y(Z,\overline{Z})=\exp{(i(Y\overline{Z}+\overline{Y}Z})+\frac{1}{2}|Y|^{2})$}
\label{1.1}
\end{equation}
Et donc :
\begin{multline}
\label{2}
\tau_{Y_{1}}\circ\tau_{Y_{2}}=\tau_{Y_{1}+Y_{2}}\exp{(-i\sigma{(Y_{1},Y_{2}}))}\ 
\\
\text{lorsqu'on identifie $Y=(y+i\eta)\in\mathbb{C}^n$ avec $Y=(y,\eta)\in\mathbb{R}^n\oplus\mathbb{R}^n$. }
\end{multline}

Soit 
\begin{equation}
\label{2.0}
f^{anti-wick}(Z,\overline{Z})=\int \exp{(i(\overline{Y}Z+\overline{Z}Y)+\frac{1}{2}|Y|^2)}m(Y)dY
\end{equation}
Alors 
\begin{equation}
F=(f^{anti-wick})^{wick}=\int \tau_Y m(Y)dY=\int \tau_Y \tilde{m}(Y)\exp{(-\frac{1}{2}|Y|^{2})}(2\pi)^{-n}dY
\end{equation}
puisqu'on posera 
\begin{equation}
\label{2.00}
\tilde{m}(Y)=\exp{(\frac{1}{2}|Y|^{2})}m(Y)(2\pi)^{n}
\end{equation}
 avec une fonction  $m(Y)\in\L^{1}(\mathbb{R}^n\oplus\mathbb{R}^n)$ , la superposition $F$ des $\tau_{Y}m(Y)dY$ donne un op\'erateur
\begin{equation}
\label{2.1}\begin{split}
&F=\int \tau_Y m(Y)dY \text{ de norme }\leq 1\\
&\text{d\`es que }\int
|m(Y)| dY\leq 1\text{ ou encore si }\| \tilde{m}\|_{L^{\infty}}\leq 1
\end{split}
\end{equation}
Composant deux op\'erateurs donn\'es par les formules (\ref{2.1}) et
relatifs \`a deux poids $m_1$ et $m_2$ int\'egrables on obtient par
composition l'op\'erateur:
\begin{equation}
\label{2.2}
F_3=F_1\circ F_2=\iint \tau_{Y_{1}+Y_{2}}\exp{(-i\sigma{(Y_{1},Y_{2}}))}m_1(Y_1)m_2(Y_2)dY_1dY_2
\end{equation}
Au moins formellement 
\begin{multline}
\label{3}
F_3=((f_3)^{anti-wick})^{wick} \\ \text{o\`u
  $f_3(Z,\overline{Z})=\iint
  \exp{(i((Y_{1}+Y_{2})\overline{Z}+\overline{{(Y_{1}+Y_{2}})}Z})+\frac{1}{2}|Y_{1}+Y_{2}|^{2})
  $} \\
\text{$\exp{((-i\sigma{(Y_{1},Y_{2}}))}m_1(Y_1)m_2(Y_2)dY_1dY_2$}
\end{multline}
On \'ecrit l'int\'egrale (\ref{3}) sous la forme : 
\begin{multline}
\label{4}
f_3(Z,\overline{Z})=\int \exp{(i(Y_{3}\overline{Z}+\overline{{Y_{3}}Z})+\frac{1}{2}|Y_{3}|^{2})}dY_3  \\ \int m_1(\frac{1}{2}Y_3+Y_4))\exp{(i\sigma ((Y_{3},Y_{4}))}m_2(\frac{1}{2}Y_3-Y_4))dY_4.
\end{multline}
avec $Y_3=Y_1+Y_2$ et $Y_4=\frac{1}{2}(Y_1-Y_2)$, si 
\begin{equation}
m_3(Y_3)=\int  m_1(\frac{1}{2}Y_3+Y_4)m_2(\frac{1}{2}Y_3-Y_4)\exp{(i\sigma{(Y_{3},Y_{4}}}))dY_4 
\label{5}
\end{equation}

Pour monter que l'int\'egrale $(\ref{5})$ converge vers vers une
fonction temp\'er\'ee en $\exp{\frac{1}{2}|Y_3|^{2}}m_3(Y_3)$ on a
besoin d'extension dans le complexe, on revient donc \`a la notation
de $Y=(y,\eta)\in \mathbb{R}^n\oplus\mathbb{R}^n$, et on fait des
prolongements dans $\mathbb{C}^n$ de $y$ et $\eta$.

On d\'eforme donc simplement le contour d'int\'egration de 

On suppose que 
\begin{itemize}
\item[i)] La fonction $\tilde{m}(y,\eta)$ se prolonge en une fonction holomorphe et born\'ee dans le domaine de $\mathbb{C}^n\oplus\mathbb{C}^n$ :
\begin{equation}
\label{6}\begin{split}
&\mathcal{D}_{\theta}=\{(y,\eta), |\Im (y,\eta)|<(\theta^2+|\Re (y,\eta)|^2)^{1/2}\} \\ &\text{pour un nombre positif $\theta$ \`a pr\'eciser}. \end{split}
\end{equation}
\end{itemize}

Soit $\Gamma_{Y_{3},\varepsilon}$ le contour d'int\'egration 
\begin{equation}
\label{7}
\tilde{y}_4=y_4+i\varepsilon \eta_3,  \tilde{\eta}_4=\eta_4-i\varepsilon y_3
\end{equation}
Les fonctions $\tilde{m}_j(\frac{1}{2}Y_3\pm Y_4))$ sont bien holomorphes dans le domaine de d\'eformation de contours si $|\varepsilon |\leq \frac{1}{2}$. 

Sur ce contour les exponentielles de  $(\ref{6})$ deviennent
oscillantes et 
on a :
\begin{multline}
\label{8}
\tilde{m}(Y_3)=\exp{(\frac{1}{4}|Y_3|^{2})}\int
\tilde{m}_1(\frac{1}{2}Y_3+(\tilde{y}_4-i\varepsilon
\eta_3,\tilde{\eta}_4+i\varepsilon
y_3))\\\tilde{m}_2(\frac{1}{2}Y_3-(\tilde{y}_4-i\varepsilon
\eta_3,\tilde{\eta}_4+ \varepsilon y_3) ) \ \\
\exp{[-(\tilde{y}_{4}^2+\tilde{\eta}_{4}^2))+\varepsilon^{2}(y_3^2+\eta_3^2)+i2\varepsilon
  (\tilde{y}_4\eta_{3}-\tilde{\eta}_4y_3)]}\\\exp{[i(\tilde{y}_4\eta_3-\tilde{\eta}_4y_3)+\varepsilon
  (\eta_{3}^{2}+y_{3}^{2}]}2^{-n}\pi^{-n}d\tilde{y}_4d\tilde{\eta}_4
\end{multline}

Le choix de $\varepsilon=-\frac{1}{2}$ est optimal et donc :

\begin{multline}
\label{9}
\tilde{m}_3(Y_3)=\\\int \tilde{m}_1(\frac{1}{2}Y_3+(\tilde{y}_4+i\frac{1}{2} \eta_3,\tilde{\eta}_4-i \frac{1}{2}y_3))\tilde{m}_2(\frac{1}{2}Y_3-(\tilde{y}_4+i\frac{1}{2} \eta_3,\tilde{\eta}_4 -i \frac{1}{2} y_3) )\\2^{-n}\pi^{-n} 
\exp{-[\tilde{y}_4^2+\tilde{\eta}_4^2]} d\tilde{y}_4d\tilde{\eta}_4
\end{multline}
 
Si on se contente d'un embryon de calcul symbolique \`a savoir que le compos\'e de deux antiwick est modulo un op\'erateur compact un op\'erateur d'antiwick on a \`a priori besoin de beaucoup moins.

On \'ecrit donc :
\begin{equation}
\label{10}
\exp{(i\sigma((Y_{3},Y_{4}))}=1+i\sigma(Y_3,Y_4)\int_0^1\exp{(i\tau \sigma((Y_{3},Y_{4}))}d\tau
\end{equation}
 Ce qui revient \`a changer dans le reste int\'egral de $(\ref{10})$; le premier terme de $(\ref{10})$ donne $F_1(Z,\overline{Z})F_2(Z,\overline{Z})$ le reste int\'egral domine un op\'erateur born\'e par :
\begin{multline}
\label{11}
\| R\|_{\mathcal{L}(F^{s},F^{s+1})}\leq\\ C\left( \int  ( | \nabla_Ym_1|+|m_1|)(Y)dY\right) \left( \int ( | \nabla_Ym_2|+|m_2|(Y))dY\right)
\end{multline}
En effet appliquer $F^{s}$ dans $F^{s+j}$ consiste \`a \'etudier les normes dans $\| \cdot \|_{\mathcal{L}(B^{s})}$ des op\'erateurs $D_{\overline{Z}^{l}}\circ B= (\overline{Z}^{\otimes_{l}}F(Z,\overline{Z}))$ qui correspondent aux fonctions $m_l(Y)=(i\partial_{Y}-\frac{\overline{Y}}{2})^{l}m(Y)$; ceci pour $0\leq l\leq j$.
\begin{corollary}
Si $N_{\alpha}=(n+|z|^2)^{\alpha})^{anti-wick}$, alors
\begin{equation}
\label{12}
N_{\alpha}\#N_{\beta}=N_{\alpha+\beta}+R \text{ $R\in \mathcal{L}(B^{s},B^{s-\Re (\alpha+\beta)-1})$}.
\end{equation}
\label{corollaire1}
\end{corollary}

On se ram\`ene en effet aux cas o\`u $\Re \alpha$ et $\Re\beta<0$ et on \'ecrit
\begin{multline}
\label{13}
N_{-\alpha}(Z,\overline{Z})=\\\iint_0^{\infty}\exp{(-(1+n)t)}t^{\alpha}\exp{i(Y\overline{Z}+\overline{Y}Z)}\frac{dt}{\Gamma (\alpha )}\exp{(-\frac{|Y|^{2}}{t})}\frac{dY}{\pi^{n} t^{n}}
\end{multline}
Le corollaire est donc \'evident.
 Cela sera g\'en\'eralement vrai pour toute superposition ad\'equates de mesures.
 
 Reste \`a traiter le cas du calcul symbolique qui implique une fonction de troncature dans dans la classe $S(1,G)$, et \`a discuter ses constantes de structure qui de toute \'evidence jouent un r\^ole crucial dans le r\'esultat.

Soit $\chi\in C^{\infty}_0(\mathbb{R}^n\oplus \mathbb{R}^n)$.

On note $\tilde{\chi}$ sa transform\'ee de Fourier pour le produit scalaire de $ \mathbb{R}^n\oplus \mathbb{R}^n$ not\'e $(X.Y)=2\Re X\overline{Y}$.

On a donc :
\begin{multline}
\label{14}
\tilde{\chi}(Y)=\int \exp{(-i(\overline{Y}Z+\overline{Z}Y))}\chi(Z)dZ\ \text{et }\chi(Z)\\=\int \exp{(i(\overline{Y}Z+\overline{Z}Y))}\tilde{\chi}(Y))\frac{dY}{(2\pi)^{2n}}
\end{multline}
Donc l'op\'erateur $\chi^{anti-wick}$ a avec les notations pr\'ec\'edentes pour fonction 
\begin{equation}
\label{15}
m(Y)=\exp{(-\frac{1}{2}|Y|^{2})}(2\pi)^{-2n}\tilde{\chi}(Y)
\end{equation}
Donc 
\begin{multline}
\label{16}
\| \chi^{anti-wck} \| \leq \int \exp{(-\frac{1}{2}|Y|^{2})}(2\pi)^{-2n}|\tilde{\chi}(Y)| dY\\\leq \sup_{Y}{|\tilde{\chi}(Y)| (2\pi)^{-n}}\leq \frac{R^{2n}\Omega_{2n-1}}{2n2^{n}\pi^{n}}
\end{multline}
Si $R=n^{1/2}\rho$, comme $\Omega_{2n-1}=\frac{\pi^{n}}{\Gamma (n)}$ on trouve que par la formule de Stirling que 
\begin{equation}
\label{17}
\| \chi^{anti-wick}\| \leq k^{n}\rho^{2n}\text{si $\| \chi\|_{\infty}\leq 1$ et $\supp \chi \subset B(0,R)$}
\end{equation}
Cette in\'egalit\'e $(\ref{17})$ montre que les normes d'op\'erateurs s'am\'eliorent avec la dimension pourvu que $\rho<k^{-1/2}$ o\`u $k$ est une constante universelle que l'on pourra calculer si n\'ecessaire.
Si $n\rightarrow \infty$ on pourra n\'egliger les op\'erateurs vivant dans des r\'egions $\supp \chi \subset B(0,\rho n^{1/2})$ et ces termes sont totalement n\'egligeables si en outre $\Lambda \leq k_2 n$ si $\rho < k_1^{-1}$.

On a donc :
\begin{proposition}
\label{prop2}
Soit $B(\Phi)=\pi(b^{anti-wick}(\Phi))$ o\`u $\pi$ est le projecteur de Bergman de $L^2(\mathbb{C}^n,\nu^{\prime})$ sur le sous-espace ferm\'e des fonctions enti\`eres de $L^2(\mathbb{C}^n,\nu^{\prime})$.

Si :
\begin{itemize}
\item[i)] $\| b^{anti-wick}\|_{L^{\infty}}\leq 1$
\item[ii)] $\supp b^{anti-wick}\subset E_{1}=\{z; |z|<\rho n^{1/2}\}$, where $\rho<\frac{1}{\sqrt{e}}$
\end{itemize}
\begin{equation}
\label{0.1}
\| B\|_{\mathcal{L}(L^{2})}\leq \| B\| _{HS}=o(n)\ \text{when $n\rightarrow \infty$}.
\end{equation}
\end{proposition}

\section{Fin de la preuve}

L'ellipticit\'e transverse entra{\^\i}ne que :
 \begin{equation}
 \label{3.1}
 a_{m}(z,\overline z)\geq C_0|z|^{m}
 \end{equation}
Il est donc naturel de prendre les symboles $S(m,G)$ avec
\begin{definition}
\label{def3}
\begin{multline}\label{2.7}
m(z)=(n^{1/2}+|z|)^{\alpha} , \ \ G_{n,z}=\frac{|dz|^{2}}{n+|z|^{2}}\ \\ \text{ et $S(m,G)$ la classe des fonctions satisfaisant \`a}\\ \|D_z^pD_{\overline z}^q f(z,\overline z)\|_{\mathcal{L}(V_{p,G},V_{q,G}^{*})}\leq C m(z)
\end{multline}
$V_{p,G}$ est $V_{p}$ muni de $G$.
\end{definition}
\begin{remark}
Les m\'etriques $G_{a,z}=\frac{|dz|^{2}}{a+|z|^{2}}$ sont bien s\^ur lentes et temp\'er\'ees d\`es que $a\geq 1$.
\end{remark}

Comme  $(\ref{3.1})$ fait que  $a^{wick}\in S((n^{1/2}+|z|)^{m},G)$; $a_{\tau}\in S((n^{1/2}+|z|)^{m-2\tau}n^{\tau},G)$, pour $0\leq \tau\leq m$. 

La condition $(\ref{3.1})$ entra{\^\i}ne que dans une zone $E_{0}=\{z; |z|\geq c_0 n^{1/2})$ o\`u $c_0$ ne d\'epend que de $a^{wick}$, on a $r_1=\sum_{\tau=1}^m a_{\tau}\in S((n^{1/2}+|z|)^{m-2},G)$. 

Soit $a^{anti-wick}=a_m(1+a_m^{-1}R_1)$, $r=a_m^{-1}r_1\in S((n^{1/2}+|z|)^{-2},G)$ dans $E_{0}$, donc si on augmente encore $c_0$,  $|r(z,\overline z)|\leq 1/2$, ce qui fait que $(1+r)^{-1}(z,\overline z)\in S(1,G)$ dans $E_{0}$. 

Dans $E_{0}$ on a bien :
\begin{equation}
\label{3.2}
a^{anti-wick}(z,\overline{z})\geq c_1 (n^{1/2}+|z|)^{m}.
\end{equation}

Les m\'etriques $G_{n,z}$ \'etant lentes elles admettent des partitions de l'unit\'e lentes \cite{Ho}. Soit $\chi_{0}\in S(1,G)$ support\'ee par $E_{0}$ on a \'egalement :
\begin{equation}
\label{3.3}
(\chi_{0} a^{anti-wick})(z,\overline{z})\geq c_1 \chi_{0} (n+|z^{2}|)^{\frac{m}{2}}
\end{equation}

Comme $N_m$ est un symbole elliptique de $S(m,G)$, on obtient aussi :
ce qui fait que 
\begin{equation}\label{3.4}
(Au,u)=((a^{anti-wick})^wu,u)\geq c_2((\chi_0n_{-\frac{m}{2}})^wu,u)-C\|u\|_{-1}^2
\end{equation}

\begin{proposition}
\label{prop2}
Soit $B(\Phi)=\pi(b^{anti-wick}(\Phi))$ o\`u $\pi$ est le projecteur de Bergman de $L^2(\mathbb{C}^n,\nu^{\prime})$ sur le sous-espace ferm\'e des fonctions enti\`eres de $L^2(\mathbb{C}^n,\nu^{\prime})$.
Si :
\begin{itemize}
\item[i)] $\| b^{anti-wick}\|_{L^{\infty}}\leq 1$
\item[ii)] $\supp b^{anti-wick}\subset E_{1}=\{z; |z|<\rho n^{1/2}\}$, where $\rho<\frac{1}{\sqrt{e}}$
\end{itemize}
\begin{equation}
\label{0.1}
\| B\|_{\mathcal{L}(L^{2})}\leq \| B\| _{HS}=o(n)\ \text{when $n\rightarrow \infty$}.
\end{equation}
\end{proposition}

\textbf{Preuve de la proposition}
Comme 
\begin{multline}
\label{0.2}
\| B\|_{HS}^2=\sum_{\alpha\in \mathbb{N}^{n}} \frac{\| B(z^{\alpha})\|^{2}}{\alpha !}\leq \int \sum_{\alpha\in\mathbb{N}^{n}}|b(z,\overline z)|^2\frac{|z|^{2\alpha}}{\alpha ! }d\nu^{\prime}(z)  \\ \leq \int_{B(0,R)}\frac{L(dz)}{\pi^{n}}=\frac{R^{2n}\Omega_{2n-1}}{2n\pi^{n}}=\frac{R^{2n}}{2n\Gamma (n)}=\frac{n^{n}\rho^{2n}}{2\Gamma{(n+1)}}
\end{multline}
La formule de Stirling prouve donc que 
\begin{equation}
\label{0.4}
\| B\|_{HS}\leq C \| b\|_{L^{\infty}}\frac{(e^{1/2}\rho)^{n}}{n^{1/4}}\ \text{$C$ \'etant ind\'ependant de $b$ et de $n$}
\end{equation}
\section{L'argument quand $n\rightarrow\infty$.}
Le calcul de Wick remplace avantageusement le calcul de Weyl pour ce qui est d'obtenir des bornes des op\'erateurs pseudo-diff\'erentiels.

Il y a essentiellement un argument qui consiste \`a dire que
l'ellipticit\'e transverse va suffire \`a assurer une borne
inf\'erieure dans une zone
$E_{I}=\{(x,\xi)\in\mathbb{R}^d\oplus\mathbb{R}^d; |(x,\xi)|\geq
C_0d^{1/2}\} $ pour une constante $C_0$ assez grande d\'ependant du
seul op\'erateur.

Le contr\^ole de $E_I^c$ est plus d\'elicat.

Il y a une troisi\`eme zone $E_3=\{(x,\xi)\in\mathbb{R}^d\oplus\mathbb{R}^d; |(x,\xi)\leq c_0d^{1/2}\}$ cette fois c'est une zone born\'ee avec pour $c_0$ une constante universelle \'egale \`a $\frac{1}{e}$, cette zone pourra aussi \^etre n\'eglig\'ee. Les difficult\'es n'arrivent que dans $E_2=\{(x,\xi)\in\mathbb{R}^d\oplus\mathbb{R}^d; c_0d^{1/2}\leq |(x,\xi)|\leq C_0d^{1/2}$, zone dans laquelle l'anti-symbole de Wick \cite{B. Lascar}, n'a pas vraiment de signe.

La variable $(x,\xi)$ \'ecrite ici est bien-sur la variable \'eclat\'ee $Y$ de la section pr\'ec\'edente, elle ne vit pour ce qui nous concernent que dans $E_C=\{(X,X_n)\in\mathbb{R}^n\oplus\mathbb{R}^n; |X|\leq C\Lambda^{1/2}\}$. Il suffit d'observer que si 
\begin{equation}
\label{16}
C\Lambda^{1/2}<\frac{1}{2}c_0d^{1/2}  \Rightarrow E_C\cap E_2=\emptyset.
\end{equation}
Ce qui sera assur\'e d\`es que $\Lambda$ et $d$ satisfont \`a :
\begin{equation}
\label{17}
\Lambda d^{-1}\leq \frac{1}{C_{1}},\  \text{avec une constante universelle $C_1$.}
\end{equation}

\begin{remerciements}
L'auteur tient \`a remercier tout d'abord le refe\'er\'ee pour ses remarques d'une aide consid\'erable ainsi que les professeurs Louis Boutet de Monvel et Johannes Sj\"ostrand pour leur aide pr\'ecieuse.
\end{remerciements}

\newpage
\appendix
\section{Comments about the paper,
{\small by R.Lascar}}

\medskip The author extends to high and
eventually infinite dimensions the results obtained with R.~Lascar in
``In\'egalit\'e de Melin H\"ormander en caract\'eristiques multiples'', to
appear in Israel Journal of Mathematics, that extended to higher
multiplicities the optimal lower bound given by L.~H\"ormander in his
book.

The statement of B.~Lascar corresponding to the ellipticity of
the Wick symbol is probably optimal.

We have not seen any errors in this paper but detailed
proofs of some statements are missing.

Our remarks are as follows:
\begin{itemize}
\item It might be a good idea to restrict the attention to the case when the
  characteristic set is a linear symplectic space on which the
  symplectic form is uniformly nondegenerate in a reasonable sense. In
  this way we have established an extension of
  the Darboux lemma to large dimensions that would allow to avoid the
  developement of a full calculus of Fourier integral operators and
  general canonical transformations.
\item The author claims that he can microlocalize--by means of a
  pseudo\-differential operator--the essential part
  of the problem to a domain in phase space,
$$E_C=\{(X,X_n);\, |X|\le C\Lambda ^{1/2},|X_n|\le C \},$$
where $\Lambda \ge 1$ is a large parameter, and the characteristic
variety $\Sigma $ is given by $X=0$, where $X=(x,\xi )\in
\mathbb{R}^{2d}$. Does there exist in the preceding works of B.~Lascar
or elsewhere in the literature the material allowing this reduction in
high dimension?
\item It seems also that one has to specify the
  assumption far from $\Sigma $. We suggest to assume
\[\begin{split}
&
\frac{1}{C}\delta ((x,\xi ),\Sigma )\le p_{m
}
(x,\xi
)\le C\delta ((x,\xi ),\Sigma ),\\
&
\delta ((x,\xi ),\Sigma )=\frac{d(x,\xi ),\Sigma }{1+d((x,\xi ),\Sigma
)},
\end{split}
\]
where $d((x,\xi ),\Sigma )$ denotes the Euclidean distance to $\Sigma
$. Far from $\Sigma $ one may then apply G\aa{}rding type arguments
(as in the recent work of J.~Nourrigat and R.~Lascar).
\item The author's remark that the domain $\{z;\, |z|\le \rho d^{1/2}
  \}$ where $\rho =\mathrm{Const.}<1/\sqrt{e}$,
is pertinent and deserves to be retained. Here $z$ is the
  complex variable corresponding to $X$ after a partial Bargmann
  transform.
\item The same is valid for the study of the region $|z|\ge
  Cd^{1/2}$,
 where the constant $C$ depends on the operator.
\item His assertions about Wick calculus, Bargman/FBI transforms, Toeplitz
  operators are valuable and interesting.

\end{itemize}


\begin{thebibliography}{5}
\bibitem{Bo}Boutet de Monvel L.\textit{ Hypoelliptic operators with double characteristics and related pseudodifferential operators}, C.P.A.M. 27 (1974) pp. 585-639.
\bibitem{Ho} H\"ormander L. \textit{The Analysis of Linear Partial Differential Operators} Volumes I-IV. Springer-Verlag. 1984-1985.
\bibitem {Lascar} Lascar B. \textit{Une classe d'op\'erateurs elliptiques du second ordre sur un espace de Hilbert.} Journal of Func. Analysis \textbf{35}, 316-343 (1980).
\bibitem{B. Lascar}Lascar B. Th\`ese de doctorat d'\'{E}tat soutenue \`a l'Universit\'e P. et M. Curie en 1978.
\bibitem{Lascars}Lascar B.- R. Lascar \textit{L'in\'egalit\'e de Melin-H\"ormander en caract\'eristiques multiples} Accept\'e pour publication. Israel Journal of maths.
\bibitem{Lerner} Lerner N. \textit{Metrics on the phase Space and Non-Selfadjoint Pseudo-Differential Operators}. Birkh\"auser Springer-Verlag. (2010).
\bibitem{Sj}Sj\"ostrand J. \textit{Singularit\'es analytiques microlocales} Ast\'erisque \textbf{95} (1984) Soci\'et\'e Math\'ematique de France. 
\end{thebibliography}
\end{document}